\documentclass [12pt] {amsart}

\usepackage [T1] {fontenc}

\usepackage {amssymb}
\usepackage {amsmath,bm}
\usepackage {amsthm}
\usepackage{tikz}
\usepackage{geometry}

\usepackage{url}

 \usepackage[titletoc,toc,title]{appendix}
 
 \usepackage{tikz-cd}

\usetikzlibrary{arrows,decorations.pathmorphing,backgrounds,positioning,fit,petri}

\usetikzlibrary{matrix}

\usepackage{enumitem}
\usepackage{hyperref}
\hypersetup{
    colorlinks=true,
    linkcolor=blue,
    filecolor=magenta,      
    urlcolor=cyan,
    citecolor=magenta
}

\DeclareMathOperator {\R} {\mathbb{R}}

\DeclareMathOperator {\Q} {\mathbb{Q}}

\DeclareMathOperator {\ws} {ws}

\DeclareMathOperator {\seq} {\subseteq}

\DeclareMathAlphabet\urwscr{U}{urwchancal}{m}{n}%
\DeclareMathAlphabet\rsfscr{U}{rsfso}{m}{n}
\DeclareMathAlphabet\euscr{U}{eus}{m}{n}
\DeclareFontEncoding{LS2}{}{}
\DeclareFontSubstitution{LS2}{stix}{m}{n}
\DeclareMathAlphabet\stixcal{LS2}{stixcal}{m} {n}

\theoremstyle {definition}
\newtheorem {definition}{Definition} [section]
\newtheorem {remark} [definition] {Remark}

\theoremstyle {plain}

\newtheorem {lemma} [definition] {Lemma}
\newtheorem {theorem} [definition] {Theorem}

\newtheorem {corollary} [definition] {Corollary}

\newtheorem {conjecture} [definition] {Conjecture}

\theoremstyle {remark}

\calclayout

\begin {document}

\title{A note on unlikely intersections in Shimura varieties}

\author{Vahagn Aslanyan and Christopher Daw}
\address{Vahagn Aslanyan, School of Mathematics, University of Leeds, Leeds, UK}

\curraddr{Department of Mathematics, University of Manchester, Manchester, UK}
\email{vahagn.aslanyan@manchester.ac.uk}

\address{Christopher Daw, Department of Mathematics and Statistics, University of Reading, Reading, UK}
\email{chris.daw@reading.ac.uk}

\thanks{{This work was supported by LMS Research in Pairs grant 42106. It was done while the first author was a senior research associate at the University of East Anglia (supported by EPSRC grant EP/S017313/1), some revisions were made while he was a research fellow at the University of Leeds (supported by Leverhulme Early Career Fellowship ECF-2022-082) and later a DKO Fellow at the University of Manchester (supported by EPSRC Open Fellowship EP/X009823/1). } }

\date{}

\vspace*{-2cm}

\begin{abstract}
We discuss the relationships between the Andr\'e-Oort, Andr\'e-Pink-Zannier, and Mordell-Lang conjectures for Shimura varieties. We then combine the latter with the geometric Zilber-Pink conjecture to obtain some new results on unlikely intersections in Shimura varieties.
\end{abstract}

\maketitle

\section{Introduction}

There are many famous conjectures and theorems in Diophantine and arithmetic geometry concerning special points on algebraic varieties. These include the Manin-Mumford and Mordell-Lang conjectures for semi-abelian varieties, and the Andr\'e-Oort and Andr\'e-Pink-Zannier conjectures for Shimura varieties. Andr\'e-Oort and Manin-Mumford are direct analogues of one another, stating, in their respective settings, that the Zariski closure of any set of special points is a finite union of special subvarieties. In \cite{Habegger-Pila-beyond}, Habegger and Pila established an analogue of the Mordell-Lang conjecture in the Shimura setting (for a product of modular curves) which was later generalised by Pila in \cite{pila-ellipt-mod-surf}. We formulate a natural generalisation of their statement for all pure Shimura varieties, and we make the observation that this generalisation is equivalent to the conjunction of the Andr\'e-Oort and Andr\'e-Pink-Zannier conjectures. This would appear timely given that proofs of these two conjectures have very recently been announced, by Pila-Shankar-Tsimerman\footnote{{The work of Pila-Shankar-Tsimerman provides the final ingredient in the full proof of Andr\'e-Oort. It was preceded by another important work due to Binyamini-Schmidt-Yafaev \cite{BSY}. In all, a large number of people have contributed to the proof in an essential way.}} \cite{Pila-Shankar-Tsimerman-andre-oort} in the case of Andr\'e-Oort, and by Richard-Yafaev \cite{Richard-Yafaev-andre-pink-2} in the case of Andr\'e-Pink-Zannier for Shimura varieties of abelian type. At the time of writing, these results are only available in preprint form. Nonetheless, we refer to them below as theorems.

\medskip

The Zilber-Pink conjecture is a vast conjecture that, when formulated for mixed Shimura varieties, encompasses all of the above. We refer the reader to \cite{Pila-ZP} for the state of the art on this conjecture. For simplicity, we work in this paper in the context of pure Shimura varieties. For the definitions, see Section \ref{sec:prelims}.

\begin{conjecture}[Zilber-Pink]\label{conj:ZP}
Let $S$ be a Shimura variety and let $V$ be a subvariety of $S$. Then $V$ contains only finitely many maximal atypical subvarieties.
\end{conjecture}

In \cite{Daw-Ren}, Ren and the second author explain a general strategy for proving Conjecture \ref{conj:ZP}, building on the earlier work of Habegger-Pila in \cite{Habegger-Pila-beyond}. The upshot is that Conjecture \ref{conj:ZP} is reduced in many cases to two arithmetic problems (sometimes referred to as the {\it large Galois orbits conjecture} and the {\it parametrisation problem}), which have since been tackled in certain cases by Orr and the second author in \cite{ExCM}, \cite{QRTUI}, \cite{PEL}, and \cite{Y(1)n}.

The unconditional result in \cite{Daw-Ren} (Proposition 6.3) is what is sometimes referred to as the {\it geometric Zilber-Pink conjecture}. (Gao refers to it as a result {\it \`a la Bogomolov}.) In this paper, as in \cite{Daw-Ren}, we combine the geometric Zilber-Pink conjecture with arithmetic to obtain a Zilber-Pink-type result. The arithmetic in this case, however, will be the aforementioned analogue of Mordell-Lang for Shimura varieties. We prove a general Zilber-Pink type result assuming the Mordell-Lang conjecture. When the latter is known, we obtain unconditional results. For instance, two of our main results are as follows. The notation and terminology are explained in Section \ref{sec:prelims}. 

\begin{theorem}\label{thm:intro-mainthm}
Let $S$ be a Shimura variety of abelian type and let $\Sigma\seq S$ denote the union of the special points of $S$ and finitely many generalised Hecke orbits. Any subvariety $V$ of $S$ contains only finitely many maximal weakly atypical subvarieties $W$ with the property that $\langle W\rangle_{\rm ws}\cap\Sigma\neq\emptyset$.
\end{theorem}

The term \textit{maximal} can be understood either as maximal among all weakly atypical subvarieties or maximal among those weakly atypical subvarieties $W$ such that $\langle W\rangle_{\rm ws}\cap\Sigma\neq\emptyset$. One can readily check that these notions coincide.

\begin{theorem}\label{thm:intro-mainthm-special}
Let $S$ be a Shimura variety and let $V$ be a subvariety of $S$. Then $V$ contains only finitely many maximal atypical subvarieties the weakly special closures of which are special subvarieties.
\end{theorem}

Observe that a weakly atypical subvariety the weakly special closure of which is special is atypical. As before, the term \textit{maximal} is unambiguous. 

For Theorem \ref{thm:intro-mainthm} we rely on the Andr\'e-Pink-Zannier conjecture for Shimura varieties of abelian type (see \cite{Richard-Yafaev-andre-pink-2}). For Theorem \ref{thm:intro-mainthm-special} we rely on the Andr\'e-Oort conjecture (see \cite{Pila-Shankar-Tsimerman-andre-oort}).

For a product of modular curves (and for semi-abelian varieties), these results were obtained by the first author in \cite{Aslanyan-remarks-atyp}. While our current proof is a generalisation of the arguments of that paper, there are some important differences, including some technicalities related to Shimura varieties. In particular, in the modular setting (i.e. in $Y(1)^n$), one can work with coordinates, which is not the case in Shimura varieties, seeing as these do not normally split as products of one-dimensional varieties. Also, we work with generalised Hecke orbits (as opposed to Hecke orbits), for the reason that the preimage of a Hecke orbit under a Shimura morphism may split into infinitely many Hecke orbits, whereas the preimage of a generalised Hecke orbit is itself a generalised Hecke orbit. This was a key observation of Richard-Yafaev (see \cite[Proposition 2.6]{Richard-Yafaev-andre-pink}). Of course, working with generalised Hecke orbits makes the result more general.

{We would like to note that two recent preprints by Barroero-Dill \cite{Bar-Dill-Hecek-ML-cat} (which is inspired by the present paper) and by Richard-Yafaev \cite{Rich-Yaf-comm-gen-AO-APZ} contain conjectures and results closely related to the ones discussed here.}

The paper is organised as follows. In Section \ref{sec:prelims} we recall some definitions from the theory of Shimura varieties and unlikely intersections. Section \ref{sec:ML} is devoted to the Andr\'e-Oort, Andr\'e-Pink-Zannier, and Mordell-Lang conjectures and the relations between them. In Section \ref{sec:ML+ZP} we show that Mordell-Lang can be combined with geometric Zilber-Pink to prove new results on unlikely intersections. In Appendix \ref{sec:appendix} we recall the relations between various notions in the theory of unlikely intersections (optimal, atypical, anomalous, etc.).

\section{Preliminaries}\label{sec:prelims}

In this section, we briefly present the necessary preliminaries on Shimura varieties. Unless stated otherwise, (sub)varieties are geometrically irreducible and identified with their sets of complex points.

\subsection{Shimura varieties}
Let $({\bf G,\bf X})$ be a Shimura datum and let $K$ be a compact open subgroup of ${\bf G}(\mathbb{A}_f)$. We assume that $\bf G$ is the generic Mumford-Tate group on $X$, seeing as this causes no loss of generality above. We obtain a complex quasi-projective algebraic variety 
\[{\rm Sh}_K({\bf G,X})={\bf G}(\Q)\backslash({\bf X}\times({\bf G}(\mathbb{A}_f)/K)),\]
which we refer to as a {\bf disconnected Shimura variety}.
We refer to any of its irreducible components as a {\bf Shimura variety} (this terminology is slightly non-standard). As such, a Shimura variety is any variety of the form $S=\Gamma\backslash X$, where $X$ is a connected component of ${\bf X}$ and $\Gamma$ is the congruence subgroup $gKg^{-1}\cap{\bf G}(\Q)_+$ for some $g\in{\bf G}(\mathbb{A}_f)$ and ${\bf G}(\Q)_+$ the subgroup of ${\bf G}(\Q)$ acting on $X$. Without loss of generality, we may always assume that $g=1$. That is, we denote by $S$ the image of $X\times\{1\}$ in ${\rm Sh}_K({\bf G},{\bf X})$.

\subsection{(Weakly) special subvarieties}\label{sec:wss}

For any Shimura subdatum $(\bf H,X_H)$ of $(\bf G,X)$ and any connected component $X_{\bf H}$ of $\bf X_H$ contained in $X$, the image of $X_{\bf H}$ in $S$ under the natural map $\pi:X\to S$ is a (closed) subvariety known as a {\bf special subvariety} of $S$. We refer to $X_{\bf H}$ as a {\bf pre-special subvariety} of $X$.

More generally, if we decompose the quotient ${\bf H}^{\rm ad}$ of $\bf H$ by its centre into a product ${\bf H}_1\times{\bf H}_2$ of two (permissibly trivial) normal $\Q$-subgroups, we obtain a decomposition $X_{\bf H}=X_1\times X_2$, and, for any $x_1\in X_1$, the image of $\{x_1\}\times X_2$ in $S$ is a subvariety known as a {\bf weakly special subvariety}. (In particular, a special subvariety is weakly special.) We refer (as in \cite[3.4]{BD}) to the decomposition $X_{\bf H}=X_1\times X_2$ as a {\bf $\Q$-splitting}.

A special subvariety of dimension zero is known as a {\bf special point}. On the other hand, any point is a weakly special subvariety. Any subvariety $W$ of $S$ is contained in a smallest weakly special subvariety of $S$, known as its {\bf weakly special closure}, denoted $\langle W\rangle_{\rm ws}$. The {\bf weakly special defect} of $W$ is defined as the difference 
\[\delta_{\rm ws}(W)=\dim\langle W\rangle_{\rm ws}-\dim W.\]

\subsection{Weakly optimal and (weakly) atypical subvarieties}\label{sec:wowas}
Let $V$ denote a subvariety of $S$. A subvariety $W$ of $V$ is known as {\bf weakly optimal} in $V$ if any subvariety of $V$ strictly containing $W$ has a strictly larger weakly special defect. In particular, $W$ is an irreducible component of $V\cap\langle W\rangle_{\rm ws}$. On the other hand, $W$ is known as {\bf (weakly) atypical} for $V$ if it is an irreducible component of the intersection of $V$ with a (weakly) special subvariety $T$ satisfying
\[\dim W>\dim V+\dim T-\dim S.\]

\subsection{Generalised Hecke orbits}
 
For any $x\in{\bf X}$ and $g\in{\bf G}(\mathbb{A}_f)$, we let $[x,g]\in {\rm Sh}_K({\bf G,X})$ denote the class of $(x,g)\in{\bf X}\times{\bf G}(\mathbb{A}_f)$. We denote by ${\bf M}={\bf MT}(x)$ the Mumford--Tate group of $x$ (that is, the smallest $\Q$-subgroup of $\bf G$ containing the image of $x:\mathbb{S}\to{\bf G}_{\R}$). By ${\rm Hom}({\bf M,G})$ we denote the set of homomorphisms defined over $\Q$. 

Following \cite[Section 2]{Richard-Yafaev-andre-pink}, we define the {\bf generalised Hecke orbit} of $x$ in ${\bf X}$ to be
\[\mathcal{H}(x)={\bf X}\cap\{\phi\circ x:\phi\in{\rm Hom}({\bf M,G})\}\]
and, for any $g\in{\bf G}(\mathbb{A}_f)$, we define the {\bf generalised Hecke orbit} of $[x,g]$ in ${\rm Sh}_K({\bf G,X})$ to be 
\[\mathcal{H}([x,g])=\{[y,h]:y\in\mathcal{H}(x),\ h\in{\bf G}(\mathbb{A}_f)\}.\]
If $[x,g]\in S$, we define its {\bf generalised Hecke orbit} in $S$ to be the intersection of $\mathcal{H}([x,g])$ with $S$.

\section{Mordell-Lang in Shimura varieties}\label{sec:ML}

Let $S$ be a Shimura variety.

\begin{definition}
 For a set $\Xi\seq S$, we denote by $\overline{\Xi}$ the union of the generalised Hecke orbits in $S$ of the points in $\Xi$. A subset $\Sigma \seq S$ is called \emph{a structure of finite rank} if there is a set $\Xi\seq S$, containing only finitely many non-special points, such that $\Sigma = \overline{\Xi}$.
\end{definition}

\begin{remark}
  { The terminology ``structure of finite rank'' is motivated by the analogy to the classical Mordell-Lang in the context of groups where its counterpart is a group of finite rank. Note that a structure of finite rank in our setting can be interpreted as a model-theoretic structure which, although not finitely generated, is sufficiently close to a finitely generated structure to be called of finite rank, though we propose no formal notion of rank in this generality.}
\end{remark}

Let $\Sigma \seq S$ be a structure of finite rank. 
 
\begin{definition}
A weakly special subvariety of $S$ is called $\Sigma$-\emph{special} if it contains a point of $\Sigma$.
\end{definition}

\begin{lemma}\label{lem:density}
A $\Sigma$-special subvariety contains a Zariski dense set of points of $\Sigma$.
\end{lemma}

\begin{proof}
Let $W$ be a $\Sigma$-special subvariety. By definition, there exists a pre-special subvariety $X_{\bf H}$ of $X$ and $\Q$-splitting $X_{\bf H}=X_1\times X_2$ (corresponding to ${\bf H}^{\rm ad}={\bf H}_1\times{\bf H}_2$) such that $W$ is equal to the image of $\{x_1\}\times X_2$ in $S$. Furthermore, there exists $x_2\in X_2$ such that the image of $x=(x_1,x_2)$ belongs to $\Sigma$. Since ${\bf H}_2(\Q)$ is dense in ${\bf H}_2(\R)$, the set $\Sigma_2={\bf H}_2(\Q)_+\cdot x_2$ is dense in $X_2$, and the image of $\{x_1\}\times \Sigma_2$ in $S$ satisfies the requirements of the lemma.
\end{proof}

The Mordell-Lang conjecture for Shimura varieties is as follows.

\begin{conjecture}[Mordell-Lang for $S$]\label{thm:mordell-lang-for-shimura}
Any subvariety $V$ of $S$ contains only finitely many maximal $\Sigma$-special subvarieties.
\end{conjecture}

Next we show that this is equivalent to the following. 

\begin{conjecture}[Mordell-Lang for $S$]\label{thm:mordell-lang-for-shimura-2}
The Zariski closure of any subset of $\Sigma$ is a finite union of $\Sigma$-special subvarieties.
\end{conjecture}

\begin{proof}[Proof that Conjecture \ref{thm:mordell-lang-for-shimura} implies Conjecture \ref{thm:mordell-lang-for-shimura-2}]
Let $\Theta\subseteq\Sigma$ and let $V=\overline{\Theta}^{\rm Zar}$. Assume that $V$ is irreducible (the general case is similar) and, without loss of generality, assume that $\Theta=V\cap\Sigma$. By Conjecture \ref{thm:mordell-lang-for-shimura}, $V$ contains finitely many maximal $\Sigma$-special subvarieties. On the other hand, since points are weakly special subvarieties, every point in $\Theta$ is contained in a maximal $\Sigma$-special subvariety of $V$.
\end{proof}

\begin{proof}[Proof that Conjecture \ref{thm:mordell-lang-for-shimura-2} implies Conjecture \ref{thm:mordell-lang-for-shimura}]
Let $V$ be a subvariety of $S$ and let $\Theta$ denote $V\cap\Sigma$. By Conjecture \ref{thm:mordell-lang-for-shimura-2}, $\overline{\Theta}^{\rm Zar}$ is equal to a finite union of $\Sigma$-special subvarieties. By Lemma \ref{lem:density}, any $\Sigma$-special subvariety contained in $V$ is contained in $\overline{\Theta}^{\rm Zar}$.
\end{proof}

Now we recall the Andr\'e-Oort and Andr\'e-Pink-Zannier conjectures. 

\begin{theorem}[Andr\'e-Oort for $S$; \cite{Pila-Shankar-Tsimerman-andre-oort}]
The Zariski closure of any set of special points is a finite union of special subvarieties. Equivalently, any subvariety $V$ of $S$ contains only finitely many maximal special subvarieties.
\end{theorem}

\begin{conjecture}[Andr\'e-Pink-Zannier for $S$]
Let $\Xi \subseteq S$ be a subset of the generalised Hecke orbit of a point in $S$. Then the Zariski closure of $\Xi$ is a finite union of $\overline{\Xi}$-special subvarieties. Equivalently, any subvariety $V$ of $S$ contains only finitely many maximal $\overline{\Xi}$-special subvarieties.
\end{conjecture}

Ordinarily the conclusion of Andr\'e-Pink-Zannier is that the irreducible components are (explicitly) only weakly special. However, since each component necessarily contains a point of $\Xi$, they are, by definition, $\overline{\Xi}$-special.

\begin{theorem}[\cite{Richard-Yafaev-andre-pink-2}]
The Andr\'e-Pink-Zannier conjecture holds when $S$ is of abelian type.
\end{theorem}

Observe that Andr\'e-Oort is a special case of Mordell-Lang, when $\Sigma$ is the set of all special points, and Andr\'e-Pink-Zannier is another special case, when $\Sigma = \overline{\{s\}}$ for an arbitrary point $s\in S$. In fact, Mordell-Lang is simply the conjunction of these two conjectures.

\begin{theorem}\label{thm:APZ->ML}
The Andr\'e-Oort and Andr\'e-Pink-Zannier conjectures imply the Mordell-Lang conjecture.
\end{theorem}

\begin{proof}
We can write $\Sigma = \Sigma_0 \cup \Sigma_1 \cup \ldots \cup \Sigma_k$, where $\Sigma_0$ contains only special points and each $\Sigma_i$ with $i=1,\ldots,k$ is equal to one generalised Hecke orbit in $S$. As such, we have 
\[ \overline{\Sigma}^{\rm{Zar}} = \bigcup_{i=0}^k  \overline{\Sigma}^{\rm{Zar}}_i. \]

By Andr\'e-Oort, $\overline{\Sigma}^{\rm{Zar}}_0$ is a finite union of special subvarieties (each of which contains a point of $\Sigma_0$). By Andr\'e-Pink-Zannier, each  $\overline{\Sigma}^{\rm{Zar}}_i$ with $i=1,\ldots,k$ is a finite union of $\Sigma$-special subvarieties. Therefore, $\overline{\Sigma}^{\rm{Zar}}$ is a finite union of $\Sigma$-special subvarieties.
\end{proof}

\begin{corollary}\label{cor:ML-for-abelian-type}
The Mordell-Lang conjecture holds when $S$ is of abelian type.
\end{corollary}

\section{Combining Mordell-Lang with geometric Zilber-Pink}\label{sec:ML+ZP}

Let $S$ be a Shimura variety, let $\Sigma \seq S$ be a structure of finite rank, and let $V$ be a subvariety of $S$.

\begin{definition}
 A weakly atypical subvariety of $V$ is called \textit{$\Sigma$-atypical} if its weakly special closure is $\Sigma$-special.
\end{definition}

\begin{lemma}\label{lem:max-atyp-is-weakly-opt}
Any maximal $\Sigma$-atypical subvariety of $V$ is weakly optimal in $V$.
\end{lemma}

\begin{proof}
Let $W \seq V$ be a $\Sigma$-atypical subvariety that is not weakly optimal. We will show that $W$ is not maximal among all $\Sigma$-atypical subvarieties of $V$. Since $W$ is not weakly optimal, there is $Y$ such that $W \subsetneq Y \seq V$ and $\delta_{\ws}(Y) \leq \delta_{\ws}(W)$. That is,
\[ \dim \langle Y \rangle_{\ws} - \dim Y \leq \dim \langle W \rangle_{\ws} - \dim W < \dim S -\dim V , \]
where the last inequality follows from the fact that $W$ is atypical. This implies that $Y$ is also an atypical subvariety of $V$. Moreover, $\langle Y \rangle_{\ws}$ contains $\langle W \rangle_{\ws}$, which contains a point of $\Sigma$. Hence, $\langle Y \rangle_{\ws}$ is $\Sigma$-special and $Y$ is $\Sigma$-atypical. We conclude that $W$ is not maximal, as required.
\end{proof}

\begin{theorem}\label{thm:ZP+ML--1}
Assume the Mordell-Lang conjecture. Every subvariety $V$ of $S$ contains only finitely many maximal $\Sigma$-atypical subvarieties.
\end{theorem}

This is equivalent to the following (the non-trivial implication can be proven exactly as in \cite[Remark 2.5]{Aslanyan-remarks-atyp}).

\begin{theorem}\label{thm:ZP+ML--2}
Assume the Mordell-Lang conjecture. For any subvariety $V$ of $S$, there is a finite collection $\Delta$ of proper $\Sigma$-special subvarieties of $S$ such that every $\Sigma$-atypical subvariety of $V$ is contained in a subvariety from $\Delta$.
\end{theorem}

\begin{proof}

Let $W$ denote a maximal $\Sigma$-atypical subvariety of $V$.

\medskip

\noindent \textbf{Step 1: $\langle W \rangle_{\ws}$ comes from finitely many families.}

\medskip
 
 By Lemma \ref{lem:max-atyp-is-weakly-opt}, $W$ is weakly optimal in $V$. By \cite[Proposition 6.3]{Daw-Ren}, there is a finite set $\Delta$ (depending only on $V$) of triples $(X_{\bf H},X_1,X_2)$, where $X_{\bf H}$ is a pre-special subvariety of $X$ and $X_{\bf H}=X_1\times X_2$ is a $\Q$-splitting, such that 
\[ \langle W \rangle_{\ws} = \pi (\{x_1\} \times X_2) \mbox{ for some } (X_{\bf H},X_1,X_2)\in\Delta\mbox{ and }x_1\in X_1. \]

\medskip

\noindent \textbf{Step 2: Pulling back to $S_{\bf H}$.} 

\medskip

Let $\Gamma_{\bf H}$ be a congruence subgroup of ${\bf H}(\Q)_+$ (the subgroup of ${\bf H}(\Q)$ acting on $X_{\bf H}$) contained in $\Gamma$, and let $\Gamma_1$ denote its image in ${\bf H}_1(\Q)$ under the natural morphism 
\[{\bf H}\to{\bf H}^{\rm ad}={\bf H}_1\times{\bf H}_2\to{\bf H_1}.\]
We obtain Shimura varieties $S_{\bf H}=\Gamma_{\bf H} \backslash X_{\bf H}$ and $S_1=\Gamma_1\backslash X_1$ and the following diagram:
\[
  \begin{tikzcd}
 S_{\bf H}  \arrow{r}{\phi} \arrow{d}{f} & S  \\
    S_1&
  \end{tikzcd}
\]

Since $\phi$ is a finite morphism, we can choose a component $\tilde{W} \seq S_{\bf H}$ of $\phi^{-1}(W)$ with $\dim \tilde{W} = \dim W$, and we choose a component $\tilde{V}$ of $\phi^{-1}(V)$ containing $\tilde{W}$. We have $\langle \tilde{W} \rangle_{\ws} = \pi_{\bf H} (\{x_1\} \times X_2)$, where $\pi_{\bf H}:X_{\bf H}\to S_{\bf H}$ denotes the natural morphism. In particular, putting $z_1=\Gamma_1x_1\in S_1$, we have $\tilde{W} \seq f^{-1}(z_1)$.

\medskip

\noindent \textbf{Step 3: There is a Zariski closed subset $V'$ of $\tilde{V}$ containing $\tilde{W}$ such that $f(V')$ is a proper Zariski-closed subset of $S_1$.} 

\medskip

Note that $f$ is closed as it extends to the (projective) Baily-Borel compactifications (see \cite{Satake}, p231). If $f(\tilde{V})\subsetneq S_1$, we choose $V'=\tilde{V}$. Suppose, therefore, that $f(\tilde{V})=S_1$. By the fibre dimension theorem (see \cite[Theorem 1.25]{Shafarevich}),
\[V'=\{ z\in \tilde{V}: \dim \tilde{V} - \dim X_1 < \dim f|_{\tilde{V}}^{-1}(f(z))\}\]
is a proper Zariski closed subset of $\tilde{V}$ and $f(V')$ is a proper closed subset of $S_1$. Since $\tilde{W}$ is weakly optimal in $\tilde{V}$, we have
\[\dim X_{\bf H}-\dim\tilde{V}\geq\dim\langle\tilde{V}\rangle_{\rm ws}-\dim\tilde{V}>\dim\langle\tilde{W}\rangle_{\rm ws}-\dim\tilde{W}=\dim X_2-\dim\tilde{W}.\]
In other words,
\[\dim\tilde{W}>\dim\tilde{V}-\dim X_{\bf H}+\dim X_2=\dim\tilde{V}-\dim X_1.\]
That is, $\tilde{W}\subseteq V'$ and we replace $V'$ with an irreducible component containing $\tilde{W}$.

\medskip

\noindent \textbf{Step 4: Applying Mordell-Lang to $f(V')$.} 

\medskip

Let $\Sigma_{\bf H} := \phi^{-1}(\Sigma)$. By \cite[Proposition 2.6]{Richard-Yafaev-andre-pink}, the preimage of a generalised Hecke orbit in $S$ is a generalised Hecke orbit in $S_{\bf H}$. Hence, the set $\Sigma_{\bf H}$ is a structure of finite rank in $S_{\bf H}$. Similarly, $\Sigma_1 = f(\Sigma_{\bf H})$ is a structure of finite rank in $S_1$. 

By Andr\'e-Pink-Zannier, $f(V')\subsetneq S_1$ contains finitely many maximal $\Sigma_1$-special subvarieties, which we denote $T_1,\ldots,T_m\subsetneq S_1$. Since $\langle \tilde{W} \rangle_{\ws}$ is $\Sigma_{\bf H}$-special, $z_1=f(\langle\tilde{W} \rangle_{\ws})$ belongs to $\Sigma_1$. Therefore, $z_1\in T_i$ for some $i\in\{1,\ldots,m\}$ and so 
\[ \tilde{W} \seq f^{-1}(z_1)=\langle \tilde{W} \rangle_{\ws} \subseteq f^{-1}(T_i), \]
from which we conclude
\[ W = \phi(\tilde{W})  \seq \phi(f^{-1}(T_i)) \subsetneq S. \]
This finishes the proof.
\end{proof}

Since Mordell-Lang is known for Shimura varieties of abelian type (see Corollary \ref{cor:ML-for-abelian-type}), the conclusion of Theorems \ref{thm:ZP+ML--1} and \ref{thm:ZP+ML--2} hold unconditionally for such Shimura varieties. This proves Theorem \ref{thm:intro-mainthm}. On the other hand, when $\Sigma$ consists of all special points, we can apply Andr\'e-Oort in Step 4 above instead of Andr\'e-Pink-Zannier (so we do not need to assume $S$ is of abelian type), yielding a proof of Theorem \ref{thm:intro-mainthm-special}.

\appendix
\section{Terminology in unlikely intersections}\label{sec:appendix}

The number of different notions (atypical, optimal, anomalous,...) appearing in the subject has become rather large over the years. For the convenience of the reader, we recall a few of these, and the relations between them. 

Let $S$ be a Shimura variety and let $V$ denote a subvariety of $S$. We have seen the notion of a weakly optimal subvariety. An {\bf optimal} subvariety of $V$ is the notion obtained by replacing the weakly special closure with the {\bf special closure} $\langle W\rangle$ (and the weakly special defect with the analogously defined {\bf special defect}).

There are also the notions of anomalous and special anomalous subvariety (see \cite{Daw-Ren}, Definitions 7.1 and 15.1): a subvariety $W$ of $V$ is {\bf anomalous} if 
\[\dim W>\max\{0,\dim V+\dim\langle W\rangle_{\ws}-\dim S\},\]
and it is {\bf special anomalous} if it satisfies the same condition with $\langle W\rangle$ in the place of $\langle W\rangle_{\rm ws}$.

The various implications between notions are laid out below. For the middle arrow on the top row, see \cite[Corollary 4.5]{Daw-Ren}. The other implications are straightforward. The label ``if pos. dim.'' means that the downward implication holds if the (weakly) atypical subvariety is positive dimensional. It has no relevance to any upward implication.

\medskip

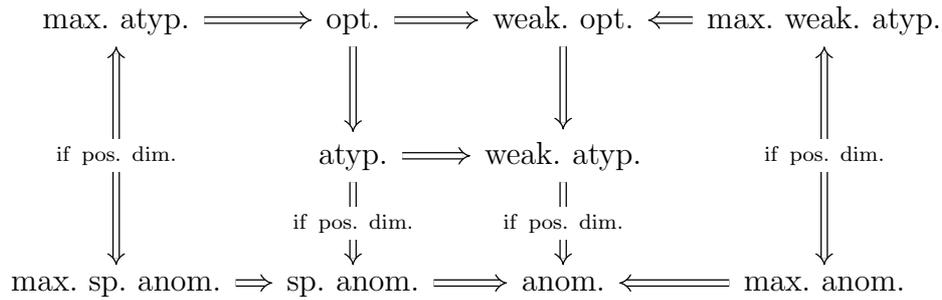
\begin{center}
\begin{tikzcd}[column sep = small, row sep = large]
{\rm max.\ atyp.} \arrow[r, Rightarrow]                                                                       & {\rm opt.} \arrow[r, Rightarrow] \arrow[d, Rightarrow]                                   & {\rm weak.\ opt.} \arrow[d, Rightarrow]                                     & {\rm max.\ weak.\ atyp.} \arrow[l, Rightarrow]                                                              \\
                                                                                                             & {\rm atyp.} \arrow[r, Rightarrow] \arrow[d, "{\rm if\ pos.\ dim.}" description, Rightarrow] & {\rm weak.\ atyp.} \arrow[d, "{\rm if\ pos.\ dim.}" description, Rightarrow] &                                                                                                           \\
{\rm max.\ sp.\ anom.} \arrow[r, Rightarrow] \arrow[uu, "{\rm if\ pos.\ dim.}" description, no head, Leftrightarrow] & {\rm sp.\ anom.} \arrow[r, Rightarrow]                                                    & {\rm anom.}                                                                & {\rm max.\ anom.} \arrow[l, Rightarrow] \arrow[uu, "{\rm if\ pos.\ dim.}" description, no head, Leftrightarrow]
\end{tikzcd}
\end{center}

\subsection*{Acknowledgement} We would like to thank the referee for some useful comments.

\bibliographystyle {alpha}
\bibliography {ref2}

\end{document}